\documentclass[a4paper,10pt]{article}

\usepackage{amssymb,amsmath,graphics,pstricks}

\newenvironment{dwd}{\par\noindent{\bf Proof.}}{\par\rightline{$\blacksquare$}}

\newtheorem{theo}{Theorem}
\newtheorem{prop}{Proposition}  
\newtheorem{coro}{Corollary}

\newtheorem{lema}{Lemma}
\newtheorem{defi}{Definition}
\def\be#1\ee{\begin{equation}#1\end{equation}}
\newcommand{\ba}{\begin{eqnarray} }
\newcommand{\ea}{\end{eqnarray} }
\def\bt#1\et{\begin{theo}#1\end{theo}}
\def\bl#1\el{\begin{lema}#1\end{lema}}
\def\bp#1\ep{\begin{prop}#1\end{prop}}
\def\bd#1\ed{\begin{defi}#1\end{defi}}
\def\Di{{\mathrm{Diam}}}
\def\ccA{{\cal A}}

\def\ccC{{\cal C}}

\def\ccE{{\cal E}}
\def\ccF{{\cal F}}

\def\ccM{{\cal M}}

\def\ccP{{\cal P}}

\def\ccS{{\cal S}}

\def\va{\varepsilon}
\def\ra{\rightarrow}

\def\E{\mathbf{E}}
\def\P{\mathbf{P}}

\def\R{{\mathbb R}}

\def\ls{\leqslant}
\def\gs{\geqslant}

\begin{document}

\title{\bf The majorizing measure approach to the sample boundedness
}
\author{Witold Bednorz
\footnote{Partially supported Research partially 
supported by MNiSW Grant no. N N201 397437.}
\footnote{Institute of Mathematics, University of Warsaw, Banacha 2, 02-097 Warszawa, Poland}}

\maketitle
\begin{abstract}
In this paper we describe the alternative approach to the sample boundedness and continuity of stochastic processes.
We show that the regularity of paths can be understood in terms of a distribution of the argument maximum.
For a centered Gaussian process $X(t)$, $t\in T$ we obtain a short proof of the exact lower bound on 
$\E \sup_{t\in T}X(t)$. Finally we prove the equivalence of a usual majorizing measure functional to its
conjugate version.  
  
\end{abstract}

\section{Introduction}\label{sect0}

Consider a Gaussian process $X(t)$, $t\in T$ on a probability space $(\Omega,\ccF,\P)$, that is a jointly Gaussian family of centered
r.v. indexed by $T$. We provide $X(t)$, $t\in T$ with the canonical distance
$$
d(s,t)=(\E(X(s)-X(t))^2)^{1/2},\;\;s,t\in T.
$$
If $X(t)$, $t\in T$ is sample bounded then 
the space $(T,d)$ has to be completely bounded since otherwise by Slepian's lemma (e.g. \cite{Le-Ta}) one can find a countable subset $S\subset T$
such that $\E\sup_{t\in S}X(t)=\infty$. It implies that $\Di (T)=\sup_{s,t\in T}d(s,t)<\infty$
and taking the Cauchy closure of $(T,d)$ one can assume that $(T,d)$ is a compact metric space. It implies that there exists
a separable modification of $X(t)$, $t\in T$ (which we refer to from now on) and therefore 
$\sup_{t\in T}X(t)$ is well defined. The sample boundedness of $X(t)$, $t\in T$
means that $\sup_{t\in T}X(t)<\infty$ almost surely. Due to the Gaussian concentration inequality the question is equivalent to the finiteness of the mean value, namely 
\be\label{nier2}
\E\sup_{t\in T}X(t)<\infty.
\ee
On the other hand note that
\be\label{emka}
\E\sup_{t\in T}X(t)=\sup_{F\subset T}\E\sup_{t\in F}X(t),
\ee
where the supremum is taken over all finite subsets $F$ of $T$. Hence (\ref{emka}) provides an alternative definition of $\E\sup_{t\in T}X(t)$, which
can be used without introducing any modification of the basic process.
\smallskip

\noindent
The second basic question on Gaussian processes is the continuity of paths. 
We say that $X(t)$, $t\in T$ is continuous if $(T,d)\ni t\ra X(t,\omega)\in \R$ is continuous, for almost all $\omega\in \Omega$.
There exists natural quantities to check whether or not the continuity takes place.  
For each $\delta>0$, define
$$
S(\delta)=\E\sup_{s,t\in T, d(s,t)\ls \delta}|X(s)-X(t)|.
$$
The if and only if condition for the continuity (see e.g. \cite{Le-Ta} Chapter 12 or \cite{Adl1} Chapter 3) is that $\lim_{\delta\ra 0}S(\delta)=0$.
\smallskip

\noindent
In this paper constant $K$ denotes a universal constant that may change from line to line.
The standard approach to the regularity of Gaussian processes goes through the entropy numbers.
Let $B(t,\va)$ be the ball of radius $\va$, centered at $t$, i.e.
$B(t,\va)=\{x\in T:\;d(x,t)\ls \va\}$. 
Denote by $N(T,d,\va)$ the smallest number of balls of radius $\va>0$ that cover $T$. The simplest upper bound of $\E \sup_{t\in T}X(t)$ was proved in \cite{Dud,Su-Tz}
$$
\E\sup_{t\in T}X(t)\ls K\int^{\infty}_0\sqrt{\log_2(N(T,d,\va))}d\va.
$$ 
Therefore $\int^{\infty}_0\sqrt{\log_2(N(T,d,\va))}d\va<\infty$
is the sufficient condition for (\ref{nier2}). It is also clear that 
$$
S(\delta)\ls K\int^{\delta}_0 \sqrt{\log_2(N(T,d,\va))}d\va,
$$
which implies the continuity of $X(t)$, $t\in T$.
Unfortunately entropy numbers does not solve the question completely, there are sample bounded Gaussian processes of
infinite entropy functional (e.g. ellipsoids in Hilbert space \cite{Tal3}) and there are discontinuous Gaussian processes that are sample bounded. 
\smallskip

\noindent
A better tool than entropies are majorizing measures. We say that a probability Borel measure $m$
is majorizing if 
\be\label{nier4}
\sup_{t\in T}\int^{\infty}_0 \sqrt{\log_2(m(B(t,\va))^{-1})}d\va<\infty.
\ee
Generalizing the notion of the majorizing measure let
$$
\ccM(\mu,\nu,\delta)=\int_T \int^{\delta}_0 \sqrt{\log_2(\mu(B(t,\va))^{-1})}d\va \nu(dt)
$$
and $\ccM(\mu,\nu)=\ccM(\mu,\nu,\Di(T))=\ccM(\mu,\nu,\infty)$.
A simple chaining argument shows (see \cite{Fer1}) that the existence of a majorizing measure
suffices for sample boundedness of $X(t)$, $t\in T$. 
\bt
The following inequality holds
$$
\E\sup_{t\in T}X(t)\ls K\inf_{\mu}\sup_{t\in T}\ccM(\mu,\delta_t), 
$$
where $\delta_t$ is the delta measure in $t$.
\et
The idea of using majorizing measures to study the sample boundedness was developed in \cite{Tal1} and \cite{Bed1}. 
In the Gaussian setting the difficult part was to prove that the existence of a majorizing measure is necessary
when $X(t)$, $t\in T$ satisfies (\ref{nier2}). The result was first proved in \cite{Tal0}.
\bt\label{theo-1} 
There exists a universal $K<\infty$ such that
$$
\E \sup_{t\in T}X(t)\gs K^{-1}\inf_{\mu}\sup_{t\in T}\ccM(\mu,\delta_t),
$$
where $K<\infty$ is a universal constant.
\et
Moreover (see e.g. \cite{Le-Ta}, Chapter 12) $X(t)$, $t\in T$ is continuous if and only if 
$$
\lim_{\delta\ra 0}\sup_{t\in T}\ccM(\mu,\delta_t,\delta)=0.
$$
A simpler argument for Theorem \ref{theo-1} appeared in \cite{Tal2}, finally in \cite{Tal3} the language of majorizing measures was replaced by admissible partitions. 
Each of the methods contains an important constructive part, where you have to 
construct a suitable admissible partition or a majorizing measure. In this paper we propose a different approach
to show the result. 
Due to \cite{Fer1} it is known that whenever $\sup_{\mu}\ccM(\mu,\mu)<\infty$ then there exists 
a majorizing measure on $T$, namely
$$
\inf_{\mu}\sup_{t\in T}\ccM(\mu,\delta_t)\ls \sup_{\mu}\ccM(\mu,\mu).
$$
The quantity $\ccM(\mu,\mu)$ is a natural upper bound for processes.
Note that for each $X(t)$, $t\in F$, where $F$
is a finite subset of $T$ there exits random $t_F$ valued in $F$ such that
$$
\E \sup_{t\in F}X(t)=\E X(t_F).
$$
Let $\mu_F(t)=\P(t_F=t)$, for $t\in F$. The measure can be treated as the distribution of the argument supremum on $F$.
We show in Section \ref{sect1} in the general setting of processes of bounded increments that $\ccM(\mu_F,\mu_F)$ is the right upper bound for the mean value of the supremum.
\bt\label{theo0}
For each $F\subset T$
$$
\E\max_{t\in F}X(t)\ls K\ccM(\mu_F,\mu_F).
$$
\et
Note that in the case of Gaussian processes the above property was proved in \cite{Adl1},
Theorem 4.2, yet also mentioned in \cite{Fer0} and known to Talagrand \cite{Tal0}.  There are many cases (see \cite{Bed2,Bed3,Bed4}) where one can prove 
the lower bound on the supremum of stochastic processes 
in the form $\sup_{\mu}\ccM(\mu,\mu)$. The benefit of the approach is that the lower bound has to be found for a given measure $\mu$ on $T$, which  
better fits the chaining argument. Moreover one can reduce the constructive part of the lower bound proof to the definition
of a natural partitioning sequence of $(T,d)$. 
Then (see Section \ref{sect4}) using this partitioning we give a short proof of the following lower bound. 
\bt\label{theo1} 
There exists a universal $K<\infty$ such that 
$$
\E\sup_{t\in F}X(t)\gs K^{-1}\sup_{\mu}\ccM(\mu,\mu).
$$
\et  
In this way we deduce that $\E\sup_{t\in T}X(t)$ is comparable with $\sup_{\mu}\ccM(\mu,\mu)$
up to a universal constant. In particular it shows the well known property 
\be\label{end0}
K^{-1}\inf_{\mu}\sup_{t\in T}\ccM(\mu,\delta_t)\ls \sup_{\mu}\ccM(\mu,\mu)\ls K\inf_{\mu}\sup_{t\in T}\ccM(\mu,\delta_t).
\ee
We prove in Section \ref{sect7} that (\ref{end0}) is true in a much generalized setting (of processes under certain increments bound).
Another question is whether or not there exists a measure $\mu_T$ such that $\E\sup_{t\in T}X(t)$ is comparable with 
$\ccM(\mu_T,\mu_T)$. Such a measure $\mu_T$ should be treated as an asymptotic argument supremum distribution, i.e. as a weak limit of
$\mu_{F_n}$ for an increasing sequence of finite $F_n$ that approximates $T$. It occurs 
that the result requires the continuity of the process.
\bt\label{theo2} 
If $X(t)$, $t\in T$ is a continuous Gaussian process then there exists measure $\mu_T$ on $T$ such that 
$$
K^{-1}\ccM(\mu_T,\mu_T)\ls \E\sup_{t\in T}X(t)\ls K\ccM(\mu_T,\mu_T).
$$
Moreover $\mu_T$ is any cluster point of any sequence $(\mu_{F_n})$ 
where $F_n\subset F_{n+1}$ and $\bigcup_n F_n$ is dense in $T$.
\et
The meaning of Theorem \ref{theo2} is that for continuous processes there exists 
an asymptotic supremum distribution, which also agrees with the result of \cite{Pol1}, where 
there is shown the existence of the argument supremum for continuous Gaussian processes at least up
to a modification of the probability space. Obviously if there exists $t_T$ such that 
$\E\sup_{t\in T}X(t)=\E X(t_T)$ then the same proof as of Theorem \ref{theo1} shows that
$\E\sup_{t\in T}X(t)\ls K\ccM(\mu_T,\mu_T)$, where $\mu_T$ is the distribution of $t_T$.
Therefore for continuous Gaussian processes there is a natural measure $\mu_T$ that 
can be used to measure $\E\sup_{t\in T}X(t)$.   
\smallskip

\noindent In the proof of Theorem \ref{theo2} we use our general estimate on $\ccS(\delta)$.
\bt\label{theo3}
There exists a universal $K<\infty$ such that
$$
K^{-1}\sup_{c>0}\sup_{\mu}(\ccM(\mu,\mu,c)-c\sqrt{\log_2(N(T,d,\delta))})\ls \ccS(\delta)\ls K\sup_{\mu}\ccM(\mu,\mu,2\delta)
$$
In particular $X(t)$, $t\in T$ is continuous if and only if $\lim_{\delta\ra \infty}\sup_{\mu}\ccM(\mu,\mu,\delta)=0$.
\et
Proof of Theorems \ref{theo2} and \ref{theo3} are provided in Section \ref{sect5}.
Then in Section \ref{sect6} we study the main toy example - Hilbert Schmidt Ellipsoid.
Finally in Section \ref{sect7} we turn to show some duality principle.
We consider the following quantity $\sup_{\mu}\inf_{t\in T}\ccM(\mu,\delta_t)$ and prove that 
in the general setting of processes of bounded increments it is comparable 
with $\inf_{\mu}\sup_{t\in T}\ccM(\mu,\delta_t)$ and hence also with $\sup_{\mu}\ccM(\mu,\mu)$.
The comparable result is discussed in a recent paper \cite{MeNa} and used to prove 
extension of the Dvoretzky theorem into the general metric spaces.
\bt\label{theo6}
There exists a universal constant $K<\infty$ such that  
$$
K^{-1}\sup_{\mu}\inf_{t\in T}\ccM(\mu,\delta_t)\ls \sup_{\mu}\ccM(\mu,\mu)\ls K\sup_{\mu}\inf_{t\in T}\ccM(\mu,\delta_t). 
$$
\et

\section{The upper bound}\label{sect1}

In this section we collect all the upper bounds required in this paper. The basic theory
was given in \cite{Tal1}  and then slightly developed in \cite{Bed1} and \cite{Bed2}.
First note that our measure approach works in much generalized setting. 
Let $(T,\rho)$ be any compact metric space and $\varphi$ - Young function, convex, increasing, 
$\varphi(0)=0$. The centered process $X(t)$, $t\in T$ is of bounded increments if
\be\label{warrunek}
\E\varphi(\frac{|X(s)-X(t)|}{\rho(s,t)})\ls 1,\;\;s,t\in T.
\ee  
Let $\Di_{\rho,}(T)$ and $B_{\rho}(t,\va)$ be diameter and ball in $\rho$ metric. 
Moreover define
$$
\sigma_{\mu,\varphi}(t,\delta)=\int^{\delta}_0 \varphi^{-1}(\frac{1}{\mu(B(t,\va))})d\va,
$$
and
$$
\ccM_{\rho,\varphi}(\mu,\nu,\delta)=\int_{T} \sigma_{\mu,\varphi}(t,\delta) \nu(dt).
$$
For simplicity let  
$$
\sigma_{\mu,\varphi}(t)=\sigma_{\mu,\varphi}(t,\Di_{\rho}(T))\;\;\mbox{and} \;\;\ccM_{\rho,\nu}(\mu,\nu)=\ccM_{\rho,\nu}(\mu,\nu,\Di_{\rho}(T)).
$$
We use the concept from the introduction i.e. let random $t_F$ valued in finite $F\subset T$ be such that 
$\E\max_{t\in F}X(t)=\E X(t_F)$ and $\mu_F(t)=\P(t_F=t)$.
\begin{prop}\label{wit4}
There exists a universal constant $K<\infty$ such that
$$
\E\sup_{t\in F}X(t) -\int_T \E X(u) \mu_F(du)\ls K\ccM_{\rho,\varphi}(\mu_F,\mu_F).
$$
\end{prop}
\begin{dwd}
First apply Theorem 1.2 from \cite{Bed1}. For each $t\in F$
the following inequality holds
\begin{align*}
&|X(t)-\int_T X(u)\mu_F(du)|\ls K_1 \sigma_{\rho,\varphi}(t)+\\
&+K_2\ccM_{\rho,\varphi}(\mu_F,\mu_F)
\int_{T\times T}\varphi(\frac{|X(u)-X(v)|}{\rho(u,v)})\nu(du,dv),
\end{align*}
where $K_1,K_2$ are absolute constants, and $\nu$ is a probability measure on $T\times T$. Denote
$$
Z=\int_{T\times T}\varphi(\frac{|X(u)-X(v)|}{\rho(u,v)})\nu(du,dv),
$$
by (\ref{warrunek}) we obtain that $\E Z\ls 1$.
Let $\Omega_t=\{t_F=t\}$, clearly
\begin{align*}
& \sum_{t\in F}\E 1_{\Omega_t}X(t)-\int_T \E X(u)\mu_F(dt)=\sum_{t\in F}\E 1_{\Omega_t}(X(t)-\int_T X(u)\mu_F(du))\ls \\
&\ls \sum_{t\in F}\E 1_{\Omega_t}(K_1\sigma_{\rho,\varphi}(t)+K_2\ccM_{\rho,\varphi}(\mu_F,\mu_F)Z)\ls K_1\sum_{t\in F}\sigma_{\rho,\varphi}(t)\mu_F(t)+\\
& +K_2\ccM(\mu_F,\mu_F)\ls  (K_1+K_2)\ccM(\mu_F,\mu_F). 
\end{align*}
It completes the proof with $K=K_1+K_2$.
\end{dwd}
We recall that in the Gaussian case i.e., when $\rho(s,t)=d(s,t)$ and $\varphi(x)=2^{x^2}-1$ we relax the notation
and use $\sigma_{\mu}$, $\ccM$ instead of $\sigma_{\rho,\varphi}$ and $\ccM_{\rho,\varphi}$.
Obviously since Gaussian variables are symmetric, $\E X(u)=0$ and hence Proposition \ref{wit4}
implies Theorem \ref{theo0}. In the non symmetric case we have the following bound. 
\begin{coro} \label{cor4}
For any $s\in T$
$$
\E\sup_{t\in T}(X(t)-X(s))\ls (1+K)\ccM(\mu_F,\mu_F).
$$
\end{coro}
\begin{dwd}
Clearly
$$
\E\sup_{t\in T}(X(t)-X(s))\ls \E\sup_{t\in T}(X(t)-\int_T X(u)\mu_F(du))+\int_{T}\E (X(u)-X(s))\mu_F(du).
$$
However $\E|X(u)-X(s)|\ls \|X(u)-X(s)\|_{\rho,\varphi}\ls  \Di_{\rho,\varphi}(T)$. The result follows 
since $\sigma_{\mu_F,\rho}(u)\gs \Di_{\rho,\varphi}(T)$.
\end{dwd}
Observe that if $\mu_{-F}$ denotes the supremum distribution of $-X(t)$ on $F$ then
\be\label{symm}
\E\sup_{t\in T}|X(t)-X(s)|\ls (1+K)(\ccM(\mu_F,\mu_F)+\ccM(\mu_{-F},\mu_{-F})).
\ee

\section{The partition structure}\label{sect2}

One of the clear consequences of Gaussian sample boundedness is that
$\Di (T)=\sup\{d(s,t):\;s,t\in T \}$ is bounded. For simplicity assume that $\Di (T)=1$. 
Recall that we apply $\sigma_{\mu}$ and $\ccM$ for this case.
\smallskip
 
\noindent 
Fix $r>1$. Let $\ccA=(\ccA_k)_{k\gs 0}$ be a partition sequence such that for each $A\in \ccA_k$
there exists $t_A\in A$ such that $A\subset B(t_A,r^{-k}/2)$.
Let $A_k(t)$ be the element of $\ccA_k$ that contains $t$.
We translate quantities $\ccM(\mu,\nu)$ into the language of $\ccA$.
\begin{lema}\label{lem2}
For each $\mu$ the following inequality holds 
\begin{align*}
& \sigma_{\mu}(t,\delta)\ls r\sum^{\infty}_{k=1}r^{-k}\sqrt{\log_2(\frac{\mu(A_{k-1}(t))}{\mu(A_k(t))})}. 
\end{align*}
\end{lema}
\begin{dwd}
First observe that
\be\label{ded1}
\int^{\infty}_0 \sqrt{\log_2(\mu(B(t,\va))^{-1})}d\va\ls  (r-1)\sum^{\infty}_{k=1}r^{-k}\sqrt{\log_2(\mu(B(t,r^{-k}))^{-1})}.
\ee
Then note that for all $t\in T$, $A_k(t)\subset B(t,r^{-k})$ and therefore
$$
\sqrt{\log_2(\mu(B(t,r^{-k}))^{-1})}\ls \sqrt{\log_2(\mu(A_k(t))^{-1})}.
$$
By the property  $\sqrt{\log_2(xy)}\ls \sqrt{\log_2(x)}+\sqrt{\log_2(y)}$ we obtain that
\be\label{wiwa1}
\sqrt{\log_2(\mu(A_k(t))^{-1})}\ls \sum^k_{l=1}\sqrt{\log_2(\frac{\mu(A_{l-1}(t))}{\mu(A_l(t))})}.
\ee
Therefore changing the summation order
\begin{align*}
&\sum^{\infty}_{k=1}r^{-k}\sqrt{\log_2(\mu(B(t,r^{-k}))^{-1})}\ls 
\sum^{\infty}_{k=1} r^{-k}\sum^{k}_{l=1}\sqrt{\log_2(\frac{\mu(A_{l-1}(t))}{\mu(A_l(t))})}=\\
&=\sum^{\infty}_{l=1} (\sum^{\infty}_{k=l}r^{-k})\sqrt{\log_2(\frac{\mu(A_{l-1}(t))}{\mu(A_l(t))})}=
\frac{r}{r-1}\sum^{\infty}_{l=1}r^{-l}\sqrt{\log_2(\frac{\mu(A_{l-1}(t))}{\mu(A_l(t))})}.
\end{align*}
It completes the proof.
\end{dwd}
\begin{coro}\label{cor2}
The following inequality holds
$$
\ccM(\mu,\nu)\ls r  \sum^{\infty}_{k=1}r^{-k}\sum_{B\in \ccA_{k-1}}\sum_{A\in \ccA_{k}(B)}\nu(A)\sqrt{\log_2(\frac{\mu(B)}{\mu(A)})}. 
$$
\end{coro}

\section{Gaussian Tools}\label{sect3}

In the general theory of Gaussian processes there are two basic properties
one can use (see Theorem 3.18 \cite{Le-Ta}, and \cite{Le} for concentration inequalities).
\bl (Sudakov Minoration)
Suppose that $d(t_i,t_j)\gs a$, $i,j\ls m$, $i\neq j$ then 
$$
\E\sup_{1\ls i\ls m}X(t_i)\gs C_1^{-1}a\sqrt{\log_2(m)},
$$
where $C_1$ is a universal constant.
\el 
\bl\label{lem1} (Gaussian Concentration)
Let $\sigma=\sup_{s,t\in D}(\E (X(t)-X(s))^2)^{1/2}$, $D\subset T$, then 
$$
\P(|\sup_{t\in D}(X(t)-\E\sup_{t\in D}X(t))|\gs u)\ls 2\exp(-\frac{u^2}{2\sigma^2}).
$$
\el
The main consequence of these facts is the basic tool we use
(see Proposition 2.1.4 in \cite{Tal3}).
\begin{prop}\label{pro}
Let $(t_i)^m_{i=1}\subset T$ satisfy $d(t_i,t_j)\gs a$ if $i\neq j$. Consider
$\sigma>0$ such that $D_i\subset B(t_i,\sigma)$. Then if $\bigcup^m_{i=1}D_i\subset D$
$$
\E \sup_{t\in D}X(t)\gs C^{-1}_1 a\sqrt{\log_2(m)}-C_2\sigma \sqrt{\log_2(m)}+\min_{1\ls i\ls m}\E\sup_{t\in D_i}X(t). 
$$ 
Thus for $a\gs (2C_1C_2)\sigma$
$$
\E\sup_{t\in D}X(t)\gs C^{-1}_3 a\sqrt{\log_2(m)}+\min_{1\ls i\ls m}\E \sup_{t\in D_i}X(t),
$$
where $C_2,C_3$ are universal constants.
\end{prop}

\section{The lower bound}\label{sect4} 

In this section we prove Theorem \ref{theo1}. Recall that $\Di(T)=1$.
First define set functionals
$$
F(A)=\E\sup_{t\in A}X(t),\;\;A\in \ccA_{k}.
$$
Using these functionals we define a natural partitioning structure for $(T,d)$. 
\smallskip

\noindent
Fix $r>1$ and $\va>0$. We construct $\ccA=(\ccA_k)_{k\gs 0}$ in the following way. 
Let $\ccA_0=\{T\}$. To define $\ccA_{k}$, $k\gs 1$ we partition each $B\in \ccA_{k-1}$ into
sets $A_1,...,A_M$ in the following way. Let $B_0=B$ and $t_1\in B$ be such that 
$$
\sup_{s\in B_0}F(C(s))\ls F(C(t_1))+\va r^{-k}, 
$$
where $C(s)= B(s,\frac{1}{2}r^{-k-1})\cap B_0$. Let $A_1=B(t_1,r^{-k}/2)$ and $B_1=B\backslash A_1$. 
We continue the construction and if
for $i\gs 1$, $B_{i-1}\neq\emptyset$ then define $t_i\in B_{i-1}$ in a way that
\be\label{imam}
\sup_{s\in B_{i-1}}F(C(s))\ls F(C(t_i))+\va r^{-k},
\ee
where $C(s)=B(s,\frac{1}{2}r^{-k-1})\cap B_{i-1}$.
Using $t_i$ we construct $A_i=B(t_i,r^{-k}/2)\cap B_{i-1}$ and $B_{i}=B_{i-1}\backslash A_i$. There exists $M<\infty$ such that
$B_M=\emptyset$, namely by the construction  $M\ls N(T,d,r^{-k}/2)<\infty$ due to the compactness of $T$.
\smallskip

\noindent
For each $B\in \ccA_k$ and $l\gs k$ denote $\ccA_l(B)=\{A\in \ccA_{l}:\;A\subset B\}$. 
Note that by the construction for each $A\in \ccA_k$ there exists $t_A\in A$
such that $A\subset B(t_A,r^{-k}/2)$, so the partition satisfies the requirement from Section \ref{sect2}.
\smallskip

\noindent
The main result of this section is the following induction scheme.
\bp\label{propo}
For $r>1$ large enough and $\va>0$ sufficiently small there exists a universal constant $L<\infty$ such that
for each measure $\mu$ on $T$ and $B\in \ccA_{k-1}$, $k\gs 1$ the following inequality holds
\begin{align*}
&\mu(B)(F(B)+4r^{-k})\gs \frac{1}{2L}r^{-k}\sum_{A\in \ccA_k(B)}\mu(A)\sqrt{\log_2(\frac{\mu(B)}{\mu(A)})} +\\
&+\sum_{C\in \ccA_{k+1}(B)}\mu(C)F(C).
\end{align*}
\ep
\begin{dwd} 
Fix $B\in\ccA_{k-1}$, $k\gs 1$. By the above construction $\ccA_k(B)=\{A_1,...,A_M\}$. 
There exits the smallest $l_0\gs 0$ such that $1\ls M\ls 2^{2^{l_0}}$, for simplicity
let $m_{-1}=0$ and $m_l=2^{2^l}$ for $l=0,1,2,...,l_0$. We group sets in $\ccA_k(B)$ using the following scheme, let $\ccA_{k,l}(B)=\{A_{m_{l-1}+1},...,A_{m_{l}}\}$,
for $0\ls l< l_0$ and $A_{k,l_0}(B)=\{A_{m_{l_0-1}+1},...,A_M\}$. Clearly $|\ccA_{k,l}(B)|=m_l-m_{l-1}$, for $0\ls l<l_0$ and 
$|\ccA_{k,l_0}(B)|= M-m_{l_0-1}$. For simplicity denote $B_l=\bigcup_{A_j\in \ccA_{k,l}(B)}A_j$, $0\ls l\ls l_0$.
\smallskip

\noindent
By the partition construction there exists points $t_i$, $1\ls i\ls M$ such that $A_i\subset B(t_i,r^{-k}/2)$ and $d(t_i,t_j)\gs r^{-k}/2$ if $1\ls i< j\ls M$.
Moreover for each $C\in \ccA_{k+1}(A_i)$ there exists $t_C\in C$ such that $C\subset B(t_C,r^{-k-1}/2)\cap A_i$ and 
hence by (\ref{imam})
$$
F(C)\ls F(D_i)+\va r^{-k}
$$
where 
$$
D_i=B(t_i,r^{-k-1}/2)\cap A_i,\;\;\mbox{for}\;1\ls i\ls M.
$$
Again by the partition construction if $i\ls j$, then 
$$
F(D_j)\ls F(D_i)+\va r^{-k}.
$$
Fix $l\gs 1$. We apply Proposition \ref{pro} with $a=r^{-k}/2$, $\sigma=r^{-k-1}/2$
and $m=m_{l-1}+1$ for sets $D_i$, $1\ls i\ls m$ and deduce that for large enough $r>1$ and 
sufficiently small $\va>0$ there exists a universal constant $L<\infty$ such that 
$$
F(B)\gs \frac{1}{L}r^{-k}2^{\frac{l}{2}}+F(C)\;\;\mbox{for all}\;C\in \ccA_{k+1}(A_j),\;A_j\in \ccA_{k,l}(B).
$$
Consequently 
\be\label{tutki1}
\mu(B_l)F(B)\gs  \frac{1}{L}\mu(B_l)r^{-k}2^{\frac{l}{2}}+\sum_{A_j\in \ccA_{k,l}(B)}
\sum_{C\in \ccA_{k+1}(A_j)}\mu(C)F(C).
\ee
The remaining bound concerns $\ccA_{k,0}(B)$. Here we cannot do better then the simplest
estimate
\be\label{tutki05}
\mu(B_0)F(B)\gs \sum_{A_j\in \ccA_{k,0}(B)}\sum_{C\in \ccA_{k+1}(A_j)}\mu(C)F(C).
\ee
By the concavity of $\sqrt{\log_2 x}$ on $[1,\infty)$ we have that for $0\ls l\ls l_0$
$$
\mu(B_l)2^{\frac{l}{2}}\gs \sum_{A_j\in A_{k,l}(B)}\mu(A_j)\sqrt{\log_2(\frac{\mu(B_l)}{\mu(A_j)})}.
$$
Moreover for each $0\ls l\ls l_0$ and $A_j\in \ccA_{k,l}(B)$
$$
\sqrt{\log_2(\frac{\mu(B_l)}{\mu(A_j)})}+\sqrt{\log_2(\frac{\mu(B)}{\mu(B_l)})}\gs \sqrt{\log_2(\frac{\mu(B)}{\mu(A_j)})}
$$
and hence
\be\label{tutki2}
\mu(B_l)(2^{\frac{l}{2}}+\sqrt{\log_2(\frac{\mu(B)}{\mu(B_l)})})\gs \sum_{A_j\in A_{k,l}(B)}\mu(A_j)\sqrt{\log_2(\frac{\mu(B)}{\mu(A_j)})}.
\ee
Thus if $2^{\frac{l}{2}}\gs \sqrt{\log_2(\frac{\mu(B)}{\mu(B_l)})}$ then by (\ref{tutki2})  
$$
2\mu(B_l)2^{\frac{l}{2}}\gs \sum_{A_j\in \ccA_{k,l}(B)}\mu(A_j)\sqrt{\log_2(\frac{\mu(B)}{\mu(A_j)})},
$$
otherwise $2^{\frac{l}{2}}\ls \sqrt{\log_2(\frac{\mu(B)}{\mu(B_l)})}$ which together with
the fact that $x\sqrt{\log_2 1+x^{-1}}$ increases on $[0,1]$ implies
$$
\mu(B_l)\sqrt{\log(\frac{\mu(B)}{\mu(B_l)})}\ls \frac{2^{\frac{l}{2}}+1}{2^{2^l}}\mu(B).
$$
Therefore due to (\ref{tutki2}) we obtain that
\be\label{tutki3}
\sum^{l_0}_{l=0}2^{\frac{l}{2}}\mu(B_l)+\sum^{l_0}_{l=0}\frac{2^{\frac{l}{2}}+1}{2^{2^l}}\mu(B)\gs 2^{-1}\sum^M_{i=1}\mu(A_i)\sqrt{\log_2(\frac{\mu(B)}{\mu(A_i)})}. 
\ee
Summing (\ref{tutki1}), (\ref{tutki05}) and (\ref{tutki3})
\begin{align*}
&\mu(B)(F(B)+r^{-k}(1+\sum^{l_0}_{l=0}\frac{2^{\frac{l}{2}}+1}{2^{2^l}}))\gs\\ 
&\gs \frac{1}{2L}r^{-k}\sum^M_{i=1}\mu(A_i)\sqrt{\log_2(\frac{\mu(B)}{\mu(A_i)})}+
\sum_{C\in \ccA_{k+1}(B)}\mu(C)F_{k+1}(C).
\end{align*}
Clearly $1+\sum^{l_0}_{l=0}\frac{2^{\frac{l}{2}}+1}{2^{2^l}}\ls 4$ which completes the proof.
\end{dwd}
Proposition \ref{propo} and the simple induction yields
$$
F(T)+4\sum^{\infty}_{k=1}r^{-2k+2}\gs \frac{1}{2L}\sum^{\infty}_{k=1}r^{-2k+1}\sum_{B\in \ccA_{2(k-1)}}\sum_{A\in \ccA_{2k-1}(B)}\mu(A)\sqrt{\log_2(\frac{\mu(B)}{\mu(A)})}.
$$
Note that for each $C\in \ccA_1$, the partition sequence $\ccA$ defines $\ccA(C)=(\ccC_k)_{k\gs 0}$ by 
$\ccA_k(C)=\ccA_{k+1}(C)$. Applying the above inequality to $C$ and $\ccA(C)$ in place of
$T$ and $\ccA$ and then using the inequality $F(T)\gs \sum_{C\in \ccA_1}\mu(C)F(C)$ we deduce that
$$
F(T)+4\sum^{\infty}_{k=1}r^{-2k+1}\gs \frac{1}{2L}\sum^{\infty}_{k=1}r^{-2k}\sum_{B\in \ccA_{2k-1}}\sum_{A\in \ccA_{2k}(B)}\mu(A)\sqrt{\log_2(\frac{\mu(B)}{\mu(A)})}.
$$
Since $F(T)=\E\sup_{t\in T}X(t)$ and for $r\gs 2$, $\sum^{\infty}_{k=1}r^{-k}\ls 1$  we finally get
$$
2(\E\sup_{t\in T}X(t)+ 4)\gs \frac{1}{2L}\sum^{\infty}_{k=1}\sum_{B\in \ccA_{k-1}}\sum_{A\in \ccA_k(B)}\mu(B)\sqrt{\log_2(\frac{\mu(B)}{\mu(A)})}.
$$
Together with Lemma \ref{lem2} and the inequality $\E\sup_{t\in T}X(t)=\E\sup_{t\in T}X(t)-X(s)\gs \sup_{t\in T}\E \max(X(t)-X(s),0)\gs C\Di (T)=C$, where 
$C$ is an absolute constant it completes the proof of Theorem \ref{theo1}.

\section{Continuity of the process}\label{sect5}

In this section we prove Theorem \ref{theo3}, i.e we show how in terms of $\sup_{\mu}\ccM(\mu,\mu,\delta)$ estimate    
$$
\ccS(\delta)=\E\sup_{s,t\in T,d(s,t)\ls \delta}|X(t)-X(s)|.
$$
For each $0<\delta\ls \Di(T)=1$ let $\ccA$ be the partition such that $A\subset B(t_A,\delta)$ for each $A\in \ccA$
and some $t_A\in A$. We require that $|\ccA|=N(T,d,\delta)$, which is clearly possible by the definition.  
Obviously
$$
\{(s,t):\:d(s,t)\ls \delta \}\supset \bigcup_{A\in \ccA} A\times \{t_A\}.
$$
Therefore we obtain that
\begin{eqnarray}\label{erth}
&& S(\delta)=\E\sup_{s,t\in T, d(s,t)\ls \delta}|X(t)-X(s)|\gs \E\max_{A\in \ccA_k}
\sup_{t\in A}|X(t)-X(t_A)|\gs \nonumber \\ 
&&\gs \sup_{A\in \ccA_k} \E\sup_{t\in A}X(t).
\end{eqnarray}
Using Theorem \ref{theo2} we get
\be\label{moon}
\E\sup_{t\in A}X(t)\gs K^{-1}\sup_{\mu_A}\ccM(\mu_A,\mu_A),
\ee
where the supremum is taken over all measures supported on $A$.
Observe that each probability measure $\mu$ on $T$ has the unique representation $\mu=\sum_{A\in \ccA}\alpha(A)\mu_A$, where $\alpha(A)\gs 1$, $\sum_{A\in \ccA}\alpha(A)=1$
and $\mu_A$ is supported on $A$. 
Consequently by the property $\sqrt{\log_2(xy)}\ls \sqrt{\log_2 x}+\sqrt{\log_2(y)}$
\begin{align*}
&\ccM(\mu,\mu,c)=\int_T\int^{c}_0 \sqrt{\log_2(\mu(B(t,\va))^{-1})} d\va\mu(dt)\ls \\
&\ls \sum_{A\in \ccA} \alpha(A)\int_T\int^{c}_0 \sqrt{\log_2(\alpha(A) \mu_A(B(t,\va)))^{-1}}\mu_A(dt)\ls \\
&\ls \sum_{A\in \ccA} [\alpha(A)\int_T\int^{\infty}_0 \sqrt{\log_2(\mu_A(B(t,\va)))^{-1}}\mu_A(dt)+c\alpha(A)\sqrt{\log_2 (\alpha(A))^{-1}}].
\end{align*}
Since by the entropy property
$$
\sum_{A\in \ccA} \alpha(A)\sqrt{\log_2(\alpha(A))^{-1}}\ls \sqrt{\log_2(N(T,d,\delta))}
$$
we deduce that 
$$
\sup_{\mu}\ccM(\mu,\mu,c)\ls \sum_{A\in \ccA}\alpha(A)\sup_{\mu_A}\ccM(\mu_A,\mu_A)+c\sqrt{\log_2(N(T,d,\delta))}. 
$$
Consequently by (\ref{erth}) and (\ref{moon})
\be\label{nock}
S(\delta)\gs K^{-1}\sup_{c>0}(\sup_{\mu}\ccM(\mu,\mu,c)-c\sqrt{\log_2(N(T,d,\delta))}).
\ee
It completes the lower bound on $S(\delta)$ in Theorem \ref{theo3}. 
\begin{coro}\label{cor5}
If $X(t)$, $t\in T$ is continuous then $\lim_{\delta\ra 0}\delta\sqrt{\log_2(N(T,d,\delta))}=0$.
\end{coro}
\begin{dwd}
Since if $X(t)$, $t\in T$ is continuous $\lim_{\delta\ra 0}S(\delta)=0$. We apply (\ref{nock})
with $c=\delta/2$ and $\mu$ equally distributed on a subset $F\subset T$ such that $|F|=N(T,d,\delta/2)$ and $\bigcup_{t\in F}B(t,\delta/2)=T$. 
By the usual entropy properties $F$ is $\delta/2$ separated and hence $\mu(B(t,\delta/2))=\mu(t)=1/(N(T,d,\delta/2))$. 
Therefore
$$
S(\delta)\gs \delta/2(\sqrt{\log_2(N(T,d,\delta/2))}-\sqrt{\log_2(N(T,d,\delta))})
$$
Using $\delta=2^{-k}$ it immediately proves that $\lim_{k\ra \infty} 2^{-k}\sqrt{\log_2(N(T,d,2^{-k}))}=0$ and hence the general result.
\end{dwd}
Consequently due to (\ref{nock}) is $X(t)$, $t\in T$ is continuous $\lim_{\delta}\sup_{\mu}\ccM(\mu,\mu,\delta)=0$.
\smallskip

\noindent
On the other hand if $F$ is the $N(T,d,\delta)$-net (i.e. $|F|=N(T,d,\delta)$, $\bigcup_{t\in F}B(t,\delta)=T$) then
$$
\E\sup_{s,t\in T}|X(t)-X(s)|\ls \E\sup_{s\in F}\sup_{t\in B(t,2\delta)}|X(t)-X(s)|.
$$
By the concentration of measure argument based on Lemma \ref{lem1} we deduce that for a universal $K_1<\infty$
\begin{align*}
&\E\sup_{s\in F}\sup_{t\in B(t,2\delta)}|X(t)-X(s)|\ls \\
&\ls \sup_{s\in F}\E\sup_{t\in B(s,2\delta)}|X(t)-X(s)|+K_1\delta\sqrt{\log_2(N(T,d,\delta))}.
\end{align*}
Since obviously
$$
\E\sup_{t\in B(s,2\delta)}|X(t)-X(s)|= 2\E\sup_{t\in B(s,2\delta)}X(t)
$$
and by Theorem \ref{theo1}, $\E\sup_{t\in B(s,2\delta)}X(t)\ls K_2\sup_{\mu}\ccM(\mu,\mu,2\delta)$ 
we deduce that
$$
\E\sup_{s\in F}\sup_{t\in B(t,2\delta)}|X(t)-X(s)|\ls K(\sup_{\mu}\ccM(\mu,\mu,2\delta)+\delta\sqrt{\log_2(N(T,d,\delta))}).
$$
Since the argument used in the proof of Corollary \ref{cor5} implies that
$$
\sup_{\mu}\ccM(\mu,\mu,2\delta)\gs \delta\sqrt{\log_2(N(T,d,\delta))},
$$
we deduce the upper bound in Theorem \ref{theo3}. Hence $\lim_{\delta}\sup_{\mu}\ccM(\mu,\mu,\delta)=0$ implies that
$\lim_{\delta\ra 0}\ccS(\delta)=0$ and therefore the process $X(t)$, $t\in T$ is continuous.
\smallskip

\noindent
We turn to prove Theorem \ref{theo2}. Assuming the continuity of $X(t)$, $t\in T$
we construct $\mu_T$ on $T$ such that 
$K^{-1}\ccM(\mu_T,\mu_T)\ls \E\sup_{t\in T}X(t)\ls K\ccM(\mu_T,\mu_T)$.
Let $(F_n)^{\infty}_{n=0}$ be any sequence of finite subsets such that
$F_n\subset F_{n+1}$ and $\bigcup_{n\gs 0}F_n$ is dense in $T$.
Due to the compactness of $(T,d)$ the sets of cluster points of $(\mu_{F_n})^{\infty}_{n=0}$ is not empty and hence 
going to a subsequence we can assume that $\mu_T$ is a weak limit of
the sequence.
Due to Theorems \ref{theo1} and \ref{theo2}
$$
K^{-1}\ccM(\mu_T,\mu_T)\ls \E\sup_{t\in T}X(t)\ls K\limsup_{n\ra \infty}\ccM(\mu_{F_n},\mu{F_n}),
$$
thus it suffices to show that
$$
\lim_{n\ra \infty} \ccM(\mu_{F_n},\mu_{F_n})=\ccM(\mu_T,\mu_T).
$$
It is clear that for $\va>0$ functionals 
$$
\Phi_{\va}(\nu)=\int_T \int^{\infty}_{\va} \sqrt{\log_2(\nu(B(t,\va))^{-1}}d\va \nu(dt)
$$
are continuous on $\ccP(T,d)$ (space of probability measures with the weak topology).
Therefore to get the convergence of $\Phi_0(\mu_{F_n})$ to $\Phi_0(\mu_T)$ we need that
$$
\sup_{n} \int_T \int^{\va}_0 \sqrt{\log_2(\mu_{F_n}(B(t,\va))^{-1})}\mu_{F_n}(dt)\ls \sup_{\mu} \ccM(\mu,\mu,\va) 
$$ 
tends to $0$, when $\va \ra 0$. Theorem \ref{theo3} implies that the above takes place whenever $X(t)$, $t\in T$ is continuous.
It completes the proof of Theorem \ref{theo2}.

\section{Hilbert Schmidt Ellipsoid}\label{sect6}

We are ready to discuss toy example for the theory of sample boundedness.
Consider $l_2$ with $\|x\|=(\sum^\infty_{i=1} x_i^2)^{1/2}$. Let  $\ccE\subset l_2$, be defined for $(t_i)^{\infty}_{i=1}$ as
$$
\ccE=\{(x_i)^{\infty}_{i=1}\in l_2:\; \sum^{\infty}_{i=1} \frac{x_i^2}{t_i^2}\ls 1\}.
$$ 
We can require that $t_i\gs t_{i+1}>0$ for $i\gs 1$. Note that $\ccE$ is compact whenever $t_i\ra 0$. Let $g=(g_i)^{\infty}_{i=1}$, where $g_i$
are independent standard Gaussian random variables. Let 
$$
X(x)=\langle x,g \rangle,\;\; x\in \ccE.
$$ 
The basic question for ellipsoid $\ccE$ is when $\sup_{x\in \ccE} X(x)<\infty$ a.s. Note that the process is continuous if sample bounded,
therefore the sample boundedness implies the existence of the supremum distribution (in the meaning of the previous section).
In the case of $\ccE$ it implies that the process $X$ may be sample bounded only if $\sum^{\infty}_{i=1}t_i^2<\infty$. Indeed 
denote for any $N<\infty$  
$$
\ccE_n=\{(x_i)^{\infty}_{i=1}:\;\sum^N_{i=1}\frac{x_i^2}{t_i^2}\ls 1\;\;\mbox{and}\;\;x_i=0\;i>N\}
$$  
Obviously $\ccE_N\subset \ccE$ and by the Schwarz inequality $\sup_{x\in \ccE_N} X(x)$ is attained on $x\in \ccE_N$ such that 
$x_i=g_i t_i^2/ (\sum^N_{i=1} t_i^2 g_i^2)^{1/2}$ for $i\ls N$ and $x_i=0$, $i>N$.
Therefore the supremum distribution $\mu_N$ on $\ccE_N$ is the distribution of $g_i t_i^2/(\sum^N_{i=1} g_i^2 t)^{1/2}$ for $i\ls N$ and $0$ for $i>N$. 
If $\sum^{\infty}_{i=1}t_i^2=\infty$ then the weak limit of $\mu_N$ is $\delta_0$ and we have a contradiction.
Therefore $\sum^{\infty}_{i=1}t_i^2<\infty$ and then the limit of $\mu_N$ exists and is equal
$\mu$  - the distribution of $g_i t_i^2/\|(g t)\|$. By Theorems \ref{theo1} and \ref{theo2} we know
that $\E \sup_{x\in \ccE} X(x)$ is comparable with $\ccM(\mu,\mu)$. 
Obviously the study of $\mu$ is a difficult question. We will provide the upper
bound on $\mu(B(x,r))$ which implies the lower bound on $\ccM(\mu,\mu)$ of the right order.
On the other hand the lower bound on $\mu(B(x,r))$ requires small value probability approach
and by now we are not able to give the right estimate. 
\smallskip

\noindent
One may note that only $\va\ls \|x\|$ are important. 
For any $y\in l^2$ define $y(i)\in l^2$ by $y(i)_j=0$, $j<i$ and $y(i)_j=y_j$ for $j\gs i$.
Denote $a_i=\|x(i)\|$, then by the construction $a_1=\|x\|$, $a_{i-1}\ls a_i$ for $i>1$ and $\lim_{i\ra\infty}a_i=0$, therefore $(a_i)^{\infty}_{i=1}$
forms a partition of $[0,\|x\|]$.  For simplicity let $a_0=t_1\gs a_1$. 
\begin{lema}
For each  $\va$ such that $a_{i+1}/\sqrt{2}\ls \va\ls a_{i}/\sqrt{2}$ for a given $i\gs 1$. 
The following inequality holds 
$$
\mu(B(x,\va))\ls \exp(-c\frac{\|t^2(i)\|^2}{t_i^4})
$$
for some universal $c>0$.
\end{lema}
\begin{dwd}
First observe that
\begin{align*}  
& \mu(B(x,\va))\ls \P(\| gt^2(i)-x(i)\|gt\| \|\ls \va \|g t\|)\ls \\
&\ls \P(\|g t^2(i)\| (\|x(i)\|^2-\va^2)^{\frac{1}{2}} \ls \langle gt^2(i),x(i)\rangle  )\ls \\
&\ls \P(\frac{1}{\sqrt{2}}\|x(i)\|\|g t^2(i)\|\ls \langle g t^2(i),x(i)\rangle ).
\end{align*}
Now there are two important tools. By Theorem 4 in \cite{Lat-Ol} for a given $\delta\in (0,1)$
and $b$ comparable with $1$ (i.e. $b_1\ls b\ls b_2$, where $b_1,b_2$ are universal constants)
\be\label{cvc1}
\P(\|g t^2(i)\|\ls \delta b\|t^2(i)\|)\ls \frac{1}{2}(2\delta)^{\frac{b\|t^2(i)\|^2}{4t^4_i}}.
\ee
The reason is that $\|gt^2(i)\|=\sup_{x\in \bar{\ccF}(i)}|\langle x(i),g\rangle|$, where $\ccF(i)=\{x\in l^2: \sum^{\infty}_{j=i}\frac{x_j^2}{t_j^4}\}$.
Therefore by the result in \cite{Lat-Ol}
$$
\P(\|gt^2(i)\|\ls \frac{1}{2}\delta M)\ls \frac{1}{2}\delta^{\frac{M^2}{4\sigma^2}}, 
$$
where $\sigma^2=\sup_{x\in \ccF(i)}\E \langle x(i), g\rangle^2$
and $M=\mathrm{Med}(\sup_{x\in \ccF(i)}\E \langle x(i), g\rangle^2)$. 
Now observe that $\sigma^2=t_i^4$ and $M$ is close to $\E \sup_{x\in \ccF(i)}|\langle x(i),g\rangle|=\E \|gt^2(i)\|$,
namely $|\E \| gt^2(i)\|-M|\ls \sigma(\pi/2)^{\frac{1}{2}}$ (e.g. \cite{Le-Ta}).
Since $\E \|gt^2(i)\|$ is comparable with $\|t^2(i)\|$ and $\sigma\ls \|t^2(i)\|$
is shows that $M$ is comparable with $\|t^2(i)\|$, i.e.
there exists $b$ comparable with $1$ such that $M=b\|t^2(i)\|$.
\smallskip

\noindent
On the other hand by the standard estimate on the gaussian measure
\begin{eqnarray}
&& \P( \frac{b\delta}{2}\|x(i)\|\|t^2(i)\|\ls \langle g t^2(i),x(i)\rangle )\ls \frac{1}{2}\exp(-\frac{b^2\delta^2}{4}\frac{\|x(i)\|^2\|t^2(i)\|^2}{\|xt^2(i)\|^2})\ls \nonumber\\
\label{cvc2}&&\ls \frac{1}{2}\exp(-\frac{b^2\delta^2}{4}\frac{\|t^2(i)\|^2}{t_i^4}). 
\end{eqnarray}
We can choose $\delta$ in a way that $\delta<1/2$ and then by (\ref{cvc1}) and (\ref{cvc2})
the lemma follows with $c=-\log(\frac{1}{2}(2\delta)^{b_1}{4}+\frac{1}{2}\exp(-\frac{b_1^2\delta^2}{2}))$,
where $b_1$ is the universal lower bound for $b$. 
\end{dwd}
Consequently $\sqrt{\log(\mu(B(x,\va))^{-1})}\gs c\frac{\|t^2(i)\|}{t_i^2}$ for
each $a_{i+1}/\sqrt{2}\ls \va\ls a_i/\sqrt{2}$. 
It suffice to prove the right lower bound on $\int_T(\|x(i)\|-\|x(i+1)\|)\mu(dx)$
\bl
The following inequality holds
$$
\int_T(\|x(i)\|-\|x(i+1)\|)\mu(dx)\gs C^{-1}\frac{t_i^4}{\|t\|\|t^2(i)\|},
$$
where $C$ is a universal constant.
\el
\begin{dwd}
First note that
$$
\|x(i)\|-\|x(i+1)\|\gs \frac{x_i^2}{2\|x(i)\|}.
$$
Now clearly
$$
\int_T \frac{x_i^2}{2\|x(i)\|}=\E \frac{t_i^4 g_i^2}{\|gt^2(i)\|\|gt\|}.
$$
Consequently
$$
\E \frac{t_i^4 g_i^2}{\|gt^2(i)\|\|gt\|}\gs \frac{t_i^4(\E |g_i|)^2}{\E\|gt^2(i)\|\|gt\|}.
$$
Then observe that $\E|g_i|=\pi$ and
$$
\E\|gt^2(i)\|\|gt\|\ls (\E \|gt^2(i)\|^2)^{1/2}(\E \|gt\|^2\|)^{1/2}=\|t^2(i)\|\|t\|.
$$
It ends the proof.
\end{dwd}
Combining Lemmas 1,2 we get
$$
\int_T\int^{a_i/\sqrt{2}}_{a_{i+1}/\sqrt{2}} \sqrt{\log(\mu(B(x,\va))^{-1})}\mu(dx)\gs \frac{1}{\sqrt{2}}cC^{-1}\frac{t_i^2}{\|t\|}.
$$
Therefore summing over $i\gs 0$ we obtain that $\ccM(\mu,\mu)\gs K^{-1}\|t\|$.
\smallskip

\noindent
Note that the idea can be extended to some class of subsets in $l^2$.

\section{The duality principle}\label{sect7}

In this section we consider general processes $X(t)$, $t\in T$ on $(T,\rho)$ under the increment condition (\ref{warrunek}).
By the result of \cite{Bed2} for $\varphi$ that satisfy 
\be\label{eg1}
\varphi(2x)\gs 2C\varphi(x)\;\;\mbox{for some}\;\; C>1 
\ee
and small enough $x\gs 0$ the following inequality holds 
$$
\ccS=\sup_{X}\E\sup_{s,t\in T}|X(t)-X(s)|=\sup_{X}\E\sup_{t\in T}X(t)\gs K^{-1}\sup_{\mu}\ccM_{\rho,\varphi}(\mu,\mu),
$$  
where the supremum is taken over all processes $X(t)$, $t\in T$ that satisfy (\ref{warrunek}).
On the other hand by (\ref{symm})
\be\label{mal2}
\ccS\ls K\sup_{\mu}\ccM_{\varphi,\rho}(\mu,\mu),
\ee
and therefore $\E\sup_{s,t\in T}|X(t)-X(s)|$ is comparable with $\sup_{\mu}\ccM_{\rho,\varphi}(\mu,\mu)$.
By the general result on majorizing measures \cite{Bed1} we have that 
$$
\ccS\ls K\inf_{\mu}\sup_{t\in T}\ccM_{\rho,\varphi}(\mu,\delta_t)
$$ 
and hence 
\be\label{mal1}
K^{-1}\inf_{\mu}\sup_{t\in T}\ccM_{\rho,\varphi}(\mu,\delta_t)\ls \sup_{\mu}\ccM_{\rho,\varphi}(\mu,\mu)\ls K\inf_{\mu}\sup_{t\in T}\ccM_{\rho,\varphi}(\mu,\delta_t),
\ee
for a large class of $\varphi$ and $\rho$ on $T$. It occurs that there is another quantity comparable with $\sup_{\mu}\ccM_{\rho,\varphi}(\mu,\mu)$, 
i.e. $\sup_{\mu}\inf_{t\in T}\ccM(\mu,\delta_t)$. We can state the main result for this section which is a generalization of Theorem \ref{theo6}.
\bt\label{theo7}
Assuming that $\varphi$ satisfies (\ref{eg1}) there exists a universal constant $K<\infty$ such that
$$
K^{-1}\sup_{\mu}\inf_{t\in T}\ccM_{\rho,\varphi}(\mu,\delta_t)\ls \sup_{\mu}\ccM_{\rho,\varphi}(\mu,\mu)\ls K\sup_{\mu}\inf_{t\in T}\ccM_{\rho,\varphi}(\mu,\delta_t)
$$
\et
\begin{dwd}
Clearly $\sup_{\mu}\ccM_{\rho,\varphi}(\mu,\mu)\gs \sup_{\mu}\inf_{t\in T}\ccM_{\rho,\varphi}(\mu,\delta_t)$,
which implies that
\be\label{mal4}
\ccS\gs K^{-1}\sup_{\mu}\inf_{t\in T}\ccM_{\varphi,\rho}(\mu,\delta_t).
\ee
Due to (\ref{mal2}) and (\ref{mal1}) for any measure $\nu$ on $T$ 
\be\label{mal3}
\ccS\ls K\sup_{t\in T}\ccM(\nu,\delta_t). 
\ee
We show that on each finite subset $F\subset T$ there exists an equality measure $\nu_F$ on $F$ such that
$\sigma_{\nu_F,\varphi}(t)$ are equal on each $t\in F$ and finite. Indeed let $F=\{t_1,...,t_m\}$,  and note that each probability measure $\mu$ on $F$ can
be treated as a point $(\alpha(1),...,\alpha(m))$ in the simplex $\triangle_m=\{(\alpha(1),...,\alpha(m)):\;\alpha(i)\gs 0, \sum^m_{i=1}\alpha(i)=1\}$, namely we set $\alpha(i)=\mu(t_i)$. 
We define mapping $\Phi:\triangle_m\ra \R^m$, $\Phi=(\Phi_1,...,\Phi_m)$ in the following way
$$
\Phi_i(\mu)=\sigma_{\mu,\varphi}(t_i)=\int^{D_{\rho,\varphi}(T)}_0 \varphi^{-1}(\frac{1}{\mu(B(t_i,\va))})d\va.
$$
\bl\label{lem5}
There exists a unique measure $\nu_F$ on $F$ such that $\Phi_i(\nu_F)$ are equal and finite for all $1\ls i\ls m$
\el
\begin{dwd}
Note that $\Phi$ is convex and continuous on $\triangle_m$. Moreover 
$\Phi_i(\delta_{t_i})=D_{\rho,\varphi}(T)$ and $\Phi_i(\delta_{t_j})=\infty$ if $i\neq j$. Therefore $\Phi$ is symplicial
in the sense that each facet of $\triangle_m$, say 
$$
\triangle_I=\{(\alpha(1),...,\alpha(n)):\;\alpha(i)=0,\; i\in I,\; \alpha(i)>0,\;i\not\in I\}
$$ 
for some $I\subset\{1,...,m\}$ is mapped on $\bar{\triangle}_I$, where
$$
\bar{\triangle}_I=\{(\Phi_1(\mu),...,\Phi_m(\mu)),\;\Phi_i(\mu)=\infty,\;i\in I,\; \Phi_i(\mu)<\infty,\;i\not\in I\}.
$$
Consequently $\triangle_{[m]}$, where $[m]=\{1,...,m\}$ must be mapped on the convex surface in $\R^m$ that connects points
$x_i=(x_i(1),...,x_m(i))$, $1\ls i \ls m$ where $x_i(i)=D_{\rho,\varphi}(T)$ and $x_i(j)=\infty$ if $i\neq j$. 
It implies that there  exists exactly one point of intersection of the surface with $y\ra (y,y,...,y)$, $y\in \R$.
Therefore there exists exactly one probability measure $\nu$ such that 
$\Phi_i(\nu)$ are equal and finite for all $1\ls i\ls m$.
\end{dwd}
Consequently by (\ref{mal3}) and Lemma \ref{lem5}
we obtain that
$$
\E\sup_{s,t\in F}|X(t)-X(s)|\ls K\inf_{t\in F}\ccM(\nu_F,\delta_t)
$$ 
and therefore 
\be\label{mal5}
\ccS\ls K\sup_{\nu}\inf_{t\in T}\ccM_{\varphi,\rho}(\nu,\delta_t).
\ee
Clearly (\ref{mal4}) and (\ref{mal5})
complete the proof.
\end{dwd}
It proves the duality principle.
\begin{coro}
The following quantities are comparable up to a universal constant: $\inf_{\mu}\sup_{t\in T}\ccM(\mu,\delta_t)$ and $\sup_{\mu}\inf_{t\in T}\ccM(\mu,\delta_t)$. 
Namely either we can search for the optimal measure $\mu$ that works for all $t\in T$ or 
for all measures we have to find the worst point $t\in T$.
\end{coro} 
As we have pointed out the result has application to the extension of the Dvoretzky 
theorem on the metric spaces.

\end{document}